\documentclass[twoside,12pt]{article}
\usepackage{amsmath,amssymb}
\begin{document}
\newtheorem{theorem}{Theorem}
\newtheorem{lemma}[theorem]{Lemma}
\newtheorem{corollary}[theorem]{Corollary}
\newtheorem{definition}[theorem]{Definition}
\newtheorem{example}[theorem]{Example}
\newtheorem{proof}{Proof}
\pagenumbering{roman}
\renewcommand{\thetheorem}{\thesection.\arabic{theorem}}
\renewcommand{\thelemma}{\thesection.\arabic{lemma}}
\renewcommand{\thedefinition}{\thesection.\arabic{definition}}
\renewcommand{\theexample}{\thesection.\arabic{example}}
\renewcommand{\theequation}{\thesection.\arabic{equation}}
\newcommand{\mysection}[1]{\section{#1}\setcounter{equation}{0}
\setcounter{theorem}{0} \setcounter{lemma}{0}
\setcounter{definition}{0}}
\newcommand{\mrm}{\mathrm}
\newcommand{\C}{\mbox{$\mathbb{C}$}}
\newcommand{\beq}{\begin{equation}}
\newcommand{\eeq}{\end{equation}}

\title
{\bf Application of a Generalized  Secant Method to Nonlinear Equations with Complex Roots}
\author
{Avram Sidi\\
Computer Science Department\\
Technion - Israel Institute of Technology\\ Haifa 32000, Israel\\
E-mail:  asidi@cs.technion.ac.il\\
http://www.cs.technion.ac.il/$\sim$asidi/}
\date{May 2021}
\maketitle
\thispagestyle{empty}
\newpage

\thispagestyle{empty}
\newpage
\pagenumbering{arabic}
\newpage
\begin{abstract}
The  secant method is a very effective numerical procedure used
for solving nonlinear equations of the form $f(x)=0$. In a recent work
[A. Sidi, Generalization of the secant method for nonlinear equations.
{\em Appl. Math. E-Notes}, 8:115--123, 2008]
 we presented a generalization of the secant method
 that uses only one  evaluation of $f(x)$ per iteration,
  and we provided a local convergence theory for it that concerns real roots.
  For each integer $k$, this method generates a sequence $\{x_n\}$ of approximations
  to a real root of $f(x)$, where, for $n\geq k$,
  $x_{n+1}=x_n-f(x_n)/p'_{n,k}(x_n)$,  $p_{n,k}(x)$ being the polynomial
  of degree $k$ that interpolates $f(x)$
  at $x_n,x_{n-1},\ldots,x_{n-k}$, the   order $s_k$ of this method satisfying $1<s_k<2$.
  Clearly, when $k=1$, this method reduces to the secant method with $s_1=(1+\sqrt{5})/2$.
  In addition,   $s_1<s_2<s_3<\cdots,$ such that
   and $\lim_{k\to\infty}s_k=2$.  In this note, we
study the application of this method to simple complex roots of a real or complex function $f(z)$.
We show that the local convergence theory developed for real roots can be extended almost as is to complex roots,
provided suitable assumptions and justifications are made.
We illustrate the theory with two numerical examples.
\end{abstract}

\vspace{2cm}

\vspace{1cm}
\noindent {\bf Keywords and expressions:} Secant method - Generalized secant method - Complex roots
\thispagestyle{empty}
\newpage
\pagenumbering{arabic}
\section{Introduction} \label{se1}
Let $\alpha$ be the solution to the equation
$f(x)=0.$
 An effective iterative method used for
solving \eqref{eq1} that makes direct use of $f(x)$ [but no
derivatives of $f(x)$] is the {\em secant method} that is
discussed in many books on numerical analysis. See, for example,
Atkinson~\cite{Atkinson:1989:INA},
Henrici~\cite{Henrici:1964:ENA}, Ralston and
Rabinowitz~\cite{Ralston:1978:FCN}, and Stoer and
Bulirsch~\cite{Stoer:2002:INA}. See also the recent note
\cite{Sidi:2006:UTR} by the author, in which the treatment of the
secant method and those of the Newton--Raphson, regula falsi, and
Steffensen methods are presented in a unified manner.

Recently, this method was generalized by the author in \cite{Sidi:2008:GSM} as follows:
Starting with $x_0,x_1,\ldots, x_k$, $k+1$ initial approximations to $\alpha$, we generate
a sequence of approximations $\{x_n\}$, via the recursion
$$ x_{n+1}=x_n-\frac{f(x_n)}{p_{n,k}'(x_n)},\quad  n=k,k+1,\ldots,$$
$p'_{n,k}(x)$ being the derivative of the polynomial $p_{n,k}(x)$  that interpolates
$f(x)$ at the points $x_n,x_{n-1},\ldots,x_{n-k}$.
Clearly, the case $k=1$  is simply the secant method. In \cite{Sidi:2008:GSM}, we also showed that
provided  $x_0,x_1,\ldots, x_k$ are sufficiently close to $\alpha$, the method
converges of order $s_k$, where $s_k\in(1,2)$ and is the only positive root of the polynomial $s^{k+1}-\sum^k_{i=0}s^i$.
We also have that
$$ \frac{1+\sqrt{5}}{2}=s_1<s_2<s_3<\cdots<2;\quad \lim_{k\to\infty}s_k=2.$$
Actually, rounded to four significant figures,
$$s_1\dot{=}1.618, \  s_2\dot{=}1.839,\ s_3\dot{=}1.928, \  s_4\dot{=}1.966, \
s_5\dot{=}1.984, \  s_6\dot{=}1.992, \  s_7\dot{=}1.996,\  \text{etc.} $$
Note that to compute  $x_{n+1}$  we need knowledge of only $f(x_n),f(x_{n-1}),\ldots,f(x_{n-k})$, and because $f(x_{n-1}),\ldots,f(x_{n-k})$  have already been  computed,  $f(x_n)$ is the only new quantity to be computed. Thus,  each step of the method requires $f(x)$ to be computed only once.
From this,  it   follows that the efficiency index of this method is simply $s_k$ and
that this index approaches $2$ by increasing $k$ even moderately.

In this work, we consider the application of this method
to simple complex roots of a function $f(z)$, where $z$ is the complex variable.
Let us denote a real or complex root of $f(z)$ by $\alpha$ again; that is, $f(\alpha)=0$ and $f'(\alpha)\neq0$.
Thus, starting with $z_0,z_1,\ldots, z_k$, $k+1$ initial approximations to $\alpha$, we generate
a sequence of approximations $\{z_n\}$ via the recursion
\beq \label{eq1} z_{n+1}=z_n-\frac{f(z_n)}{p_{n,k}'(z_n)},\quad  n=k,k+1,\ldots,\eeq
$p_{n,k}'(z)$ being the derivative of the polynomial $p_{n,k}(z)$ that
interpolates $f(z)$ at the points $z_n,z_{n-1},\ldots,z_{n-k}$. As in \cite{Sidi:2008:GSM}, we can use
Newton's interpolation  formula to generate
$p_{n,k}(z)$ and $p_{n,k}'(z)$.
Thus
\beq
\label{eq2}p_{n,k}(z)=f(z_n)+\sum^{k}_{i=1}f[z_n,z_{n-1},\ldots,z_{n-i}]
\prod^{i-1}_{j=0}(z-z_{n-j}).\eeq
and
\beq
\label{eq3} p_{n,k}'(z)=f[z_n,z_{n-1}]+\sum^{k}_{i=2}f[z_n,z_{n-1},\ldots,z_{n-i}]
\prod^{i-1}_{j=1}(z_n-z_{n-j}).\eeq
Here $g[\zeta_0,\zeta_1,\ldots,\zeta_m]$ is the divided difference of order $m$ of the function $g(z)$ over the set of points $\{\zeta_0,\zeta_1,\ldots,\zeta_m\}$ and is a symmetric function of these points.
For details, we refer the reader to \cite{Sidi:2008:GSM}.

As proposed in \cite{Sidi:2008:GSM}, we generate the $k+1$ initial approximations as follows: We choose the  approximations $z_0,z_1$ first. We then generate $z_2$ by applying our method with $k=1$ (that is, with the secant method). Next, we apply our method to $z_0,z_1,z_2$  with $k=2$ and obtain $z_3$, and so on, until we have generated all $k$ initial approximations, via
$$ z_{j+1}=z_j-\frac{f(z_j)}{p_{j,j}'(z_j)}, \quad j=1,2,\ldots,k-1.$$

\noindent{\bf Remarks.}
\begin{enumerate}
\item
 Instead of {\em choosing}  $z_1$ arbitrarily, we can {\em generate} it as $z_1=z_0+f(z_0)$ as suggested in Brin \cite{Brin:2016:TTNA}, which is quite sensible since $f(z)$ is small near the root $\alpha$. We can also use the method of Steffensen---which uses only $f(z)$ and no derivatives of $f(z)$---to generate $z_1$ from $z_0$; thus,
 $$ z_1=z_0-\frac{[f(z_0)]^2}{f(z_0+f(z_0))-f(z_0)}.$$

\item
It is clear that, in case $f(z)$ takes on only real values along the $\text{Re}\,z$ axis
 and we are looking for nonreal roots of $f(z)$,
at least one of the initial approximations must be chosen to be nonreal.
\end{enumerate}

In the next section, we analyze the local convergence properties of the method as
 it is applied to complex roots. We show that the analysis of \cite{Sidi:2008:GSM}
can be extended to the complex case following some clever manipulation.
 We prove that the order $s_k$ of the method
 is the same as that we discovered in the real case.
In Section \ref{se3}, we provide two numerical examples to confirm the results
of our convergence analysis.

\section{Local convergence analysis}\label{se2}
\setcounter{theorem}{0} \setcounter{equation}{0}
 We now turn to
the analysis of the sequence $\{z_n\}^{\infty}_{n=0}$ that is
generated via \eqref{eq1}. We treat the
case $k\geq 2$. The case  $k=1$ is similar but much simpler.

In our analysis, we will make use of the Hermite-Genocchi formula that provides an integral
representation for divided differences.\footnote{For a proof of this formula,
see Atkinson \cite{Atkinson:1989:INA}, for example.} Even though this formula is usually stated for
 functions defined on real intervals,
it is easy to verify (see Filipsson \cite{Filipsson:1997:CMV}, for example) that it also applies to functions defined in the complex plane under proper assumptions. Thus,
provided $g\in C^m(E)$, $E$ being a bounded closed convex set in the complex plane, and provided
$\zeta_0,\zeta_1,\ldots, \zeta_m$ are in $E$, there holds
\beq \label{eqHG} g[\zeta_0,\zeta_1,\ldots, \zeta_m]=\idotsint\limits_{S_m}
                               g^{(m)}(t_0\zeta_0+t_1\zeta_1+\cdots+t_m\zeta_m)\,dt_1\cdots dt_m.\eeq
Here $S_m$ is the $m$-dimensional simplex defined as
   \beq   \label{eqSm}   S_m = \bigg\{(t_1,\ldots,t_m):\ t_i\geq0,\ i=1,\ldots,m,\quad
                   \sum_{i=1}^m t_i\leq1\bigg\}, \eeq
and
\beq \label{eqt0} t_0=1-\sum_{i=1}^m t_i.\eeq
We also note that \eqref{eqHG} holds whether  the $\zeta_i$ are  distinct or not.

   By the conditions we have imposed on $g(z)$ it is easy to see that the integrand
$g^{(m)}(\sum^m_{i=0}t_i\zeta_i)$ in \eqref{eqHG} is always defined because  $\sum^m_{i=0}t_i\zeta_i$ is in the set $E$ and $g\in C^m(E)$. This is so because, by \eqref{eqSm} and \eqref{eqt0},
$$ (t_1,\ldots,t_m)\in S_m \quad \Rightarrow \quad t_i\geq0,\quad i=0,1,\ldots,m, \quad \text{and}\quad
\sum^m_{i=0}t_i=1,$$  which implies that  $\sum^m_{i=0}t_i\zeta_i$ is a convex combination of
$\zeta_0,\zeta_1,\ldots, \zeta_m$ hence is
 in the set $C=\text{conv}\{\zeta_0,\zeta_1,\ldots, \zeta_m\}$, the convex hull of the points $\zeta_0,\zeta_1,\ldots, \zeta_m$,  and $C\subset E$.
  Consequently, taking moduli on both sides of \eqref{eqHG}, we obtain, for all $\zeta_i$ in $E$,
\begin{align}\label{eqHGb}\big |g[\zeta_0,\zeta_1,\ldots, \zeta_m]\big|&\leq
\idotsint\limits_{S_m}\bigg| g^{(m)}\bigg(\sum^m_{i=0}t_i\zeta_i\bigg)\bigg|\,dt_1\cdots dt_m   \notag \\
 &\leq\frac{\|g^{(m)}\|}{m!},\quad
\|g^{(m)}\|=\max_{z\in E}|g^{(m)}(z)|.\end{align}
In addition, since $\sum^m_{i=0}t_i=1$ in \eqref{eqHG},
as $\zeta_i\to\hat{\zeta}$, for all $i=0,1,\ldots,m$, there hold
$\sum^m_{i=0}t_i\zeta_i\to \hat{\zeta}$ and
$g^{(m)}(\sum^m_{i=0}t_i\zeta_i)\to g^{(m)}(\hat{\zeta})$, and hence
\beq \label{eqlim} \lim_{\substack{\zeta_i\to \hat{\zeta}\\ i=0,1,\ldots,m}}g[\zeta_0,\zeta_1,\ldots, \zeta_m]=g[\underbrace{\hat{\zeta},\hat{\zeta},\ldots,\hat{\zeta}}_{m+1\ \text{times}}]=\frac{g^{(m)}(\hat{\zeta})}{m!}.\eeq
In \eqref{eqHGb} and \eqref{eqlim},  we have also invoked the fact that (see \cite[p. 346]{Davis:1984:MNI}, for example)
$$ \idotsint\limits_{S_m} dt_1\cdots dt_m=\frac{1}{m!}.$$
We will make use of these in the proof of our main theorem that follows.
This theorem and its proof are almost identical to that given in \cite{Sidi:2008:GSM},
but taking into account where and when needed, the fact that we are now working in the complex plane.
For convenience, we provide all the details of the proof.

\begin{theorem}\label{th1}
Let $\alpha$ be a simple root of $f(z)$, that is, $f(\alpha)=0$, but $f'(\alpha)\neq0$.
Assume that
$f\in C^{k+1}(B_r)$, where $B_r$  is a closed ball of radius $r$ containing
$\alpha$ as its center, that is,
\beq B_r=\{z\in \C:\ |z-\alpha|\leq r\}.\eeq
Let $z_0,z_1,\ldots,z_k$ be distinct initial
approximations to $\alpha$, and generate $z_n$,
$n\geq k+1$,  via
 \begin{equation} \label{eq:gs}
z_{n+1}=z_n-\frac{f(z_n)}{p'_{n,k}(z_n)},\quad
n=k,k+1,\ldots,\end{equation} where $p_{n,k}(z)$ is the polynomial
of interpolation to $f(z)$ at the points
$z_n,z_{n-1},\ldots,z_{n-k}$.
Then, provided $z_0,z_1,\ldots,z_k$ are in $B_r$ and sufficiently close to
$\alpha$, we have the following  cases:
\begin{enumerate}
\item
 If $f^{(k+1)}(\alpha)\neq0$, the sequence $\{z_n\}$ converges to $\alpha$, and

\begin{equation}\label{eq:130}
\lim_{n\to\infty}\frac{\epsilon_{n+1}}{\prod^k_{i=0}\epsilon_{n-i}}=
\frac{(-1)^{k+1}}{(k+1)!}\,\frac{f^{(k+1)}(\alpha)}{f'(\alpha)}\equiv
L;\quad \epsilon_n=z_n-\alpha\quad \forall\, n. \end{equation}
 The order of convergence is $s_k$, $1<s_k<2$, where $s_k$ is the only
 positive root of the equation $s^{k+1}=\sum_{i=0}^k s^i$ and
 satisfies
\begin{equation}\label{eq:132} 2-2^{-k-1}e< s_k<2-2^{-k-1}\quad\text{for $k\ge 2$;} \quad s_k<s_{k+1};\quad
\lim_{k\to\infty}s_k=2,\end{equation}   $e$ being  the base of
natural logarithms, and
\begin{equation}\label{eq:135}
\lim_{n\to\infty}\frac{|\epsilon_{n+1}|}{~~~|\epsilon_{n}|^{s_k}}=|L|^{(s_k-1)/k}\quad\Rightarrow\quad s_k=\lim_{n\to\infty}\frac{\log|\epsilon_{n+1}/\epsilon_{n}|}
{\log|\epsilon_{n}/\epsilon_{n-1}|}.\end{equation}
\item
 If $f(z)$ is a polynomial of degree at most $k$, the sequence $\{z_n\}$ converges to $\alpha$, and
\begin{equation}\label{eq:177}
\lim_{n\to\infty}\frac{\epsilon_{n+1}}{\epsilon_n^2}=
\frac{f''(\alpha)}{2f'(\alpha)};\quad \epsilon_n=z_n-\alpha\quad \forall\ n. \end{equation}
Thus $\{z_n\}$ converges  of order $2$ if $f''(\alpha)\neq0$, and of  order  greater than $2$ if $f''(\alpha)=0$.

\end{enumerate}
\end{theorem}

\noindent{\bf Proof.}
We start by deriving a closed-form expression for the error in
$z_{n+1}$.  Subtracting $\alpha$ from both sides of \eqref{eq:gs},
and noting that
$$ f(z_n)=f(z_n)-f(\alpha)=f[z_n,\alpha](z_n-\alpha),$$
we have
\begin{equation}\label{eq:112}
z_{n+1}-\alpha=\bigg(1-\frac{f[z_n,\alpha]}{p'_{n,k}(z_n)}\bigg)(z_n-\alpha)=
\frac{p'_{n,k}(z_n)-f[z_n,\alpha]}{p'_{n,k}(z_n)}\,(z_n-\alpha).
\end{equation}
 We now note that
\beq \label{eq456}
p'_{n,k}(z_n)-f[z_n,\alpha]=\big\{p'_{n,k}(z_n)-f'(z_n)\big\}+
\big\{f'(z_n)-f[z_n,\alpha]\big\},\eeq  and that
\beq\label{eq:107}
f'(z_n)-p'_{n,k}(z_n)=f[z_n,z_n,z_{n-1},\ldots,z_{n-k}]\prod^k_{i=1}(z_n-z_{n-i})\eeq
and
\beq \label{eq:107a}
f'(z_n)-f[z_n,\alpha]=f[z_n,z_n]-f[z_n,\alpha]=f[z_n,z_n,\alpha](z_n-\alpha).\eeq
For simplicity of notation, let
\beq\label{eq123} -f[z_n,z_n,z_{n-1},\ldots,z_{n-k}]=\widehat{D}_n\quad\text{and}\quad f[z_n,z_n,\alpha]=\widehat{E}_n, \eeq
and rewrite \eqref{eq456} and \eqref{eq:107}  as
\begin{align} &p'_{n,k}(z_n)-f[z_n,\alpha]=\widehat{D}_n
\prod^k_{i=1}(\epsilon_n-\epsilon_{n-i})+\widehat{E}_n\epsilon_n,\\
&p'_{n,k}(z_n)=f'(z_n)+\widehat{D}_n\prod^k_{i=1}(\epsilon_n-\epsilon_{n-i}).
\end{align}
Substituting  these in
\eqref{eq:112},  we finally
obtain

\begin{equation}\label{eq:123}
\epsilon_{n+1}=C_n\epsilon_n;\quad
C_n\equiv\frac{p'_{n,k}(z_n)-f[z_n,\alpha]}{p'_{n,k}(z_n)}=\frac{\widehat{D}_n\prod^k_{i=1}(\epsilon_n-\epsilon_{n-i})
+\widehat{E}_n
\epsilon_n}{f'(z_n)+\widehat{D}_n\prod^k_{i=1}(\epsilon_n-\epsilon_{n-i})}.\end{equation}

We now prove that convergence takes place. First, let us  assume without loss of generality that $f'(z)\neq0$ for all $z\in B_r$, and set $m_1=\min_{z\in B_r}|f'(z)|>0$.
[This is
possible since $\alpha\in B_{r}$ and $f'(\alpha)\neq 0,$ and we can choose $r$ as small as we wish
to also guarantee $m_1>0$.]
Next, let $M_s=\max_{z\in B_r}|f^{(s)}(z)|/s!,$ $s=1,2,\ldots.$
Thus, assuming that $\{z_n,z_{n-1},\ldots,z_{n-k}\}\subset B_r$, we have  by \eqref{eqHGb} that
$$ |\widehat{D}_n|\leq M_{k+1}, \quad |\widehat{E}_n|\leq M_2,\quad \text{because $\{\alpha,z_n,z_{n-1},\ldots,z_{n-k}\}\subset B_r$}.$$
Next,  choose the  ball $B_{t/2}$ sufficiently small (with $t/2\leq r$) to ensure that
$m_1>2M_{k+1}t^k+M_2t/2$.  It can now be
shown that, provided $z_{n-i}$, $i=0,1,\ldots,k,$ are all in $B_{t/2}$,
there holds

\begin{align*}|C_n|&\leq
\frac{M_{k+1}\prod^k_{i=1}|\epsilon_n-\epsilon_{n-i}|+M_2|\epsilon_n|}
{m_1-M_{k+1}\prod^k_{i=1}|\epsilon_n-\epsilon_{n-i}|}  \\
&\leq
\frac{M_{k+1}\prod^k_{i=1}(|\epsilon_n|+|\epsilon_{n-i})|+M_2|\epsilon_n|}
{m_1-M_{k+1}\prod^k_{i=1}(|\epsilon_n|+|\epsilon_{n-i}|)}\leq
\overline{C},\end{align*} where
$$
\overline{C}\equiv\frac{M_{k+1}t^k+M_2t/2}{m_1-M_{k+1}t^k}<1.$$
Consequently, by \eqref{eq:123}, $|\epsilon_{n+1}|<\overline{C}|\epsilon_n|<|\epsilon_n|$,
which implies that $z_{n+1}\in B_{t/2}$, just like $z_{n-i}$,
$i=0,1,\ldots,k.$ Therefore, if  $z_0,z_1,\ldots,z_k$ are chosen
in $B_{t/2}$, then $|C_n|\leq \overline{C}<1$ for all $n\geq k$, hence
$\{z_n\}\subset B_{t/2}$  and $\lim_{n\to\infty}z_n=\alpha$.

As for \eqref{eq:130} when $f^{(k+1)}(\alpha)\neq0$, we proceed as follows: By the fact that
$\lim_{n\to\infty}z_n=\alpha$,  we first note that, by \eqref{eq:107} and \eqref{eq:107a},

\beq \label{ppff} \lim_{n\to\infty}p'_{n,k}(z_n)=f'(\alpha)=\lim_{n\to\infty}f[z_n,\alpha],\eeq
and thus $\lim_{n\to\infty}C_n=0$. This means that
$\lim_{n\to\infty}(\epsilon_{n+1}/\epsilon_n)=0$ and,
equivalently, that $\{z_n\}$ converges of order greater than 1.
  As a result,

$$\lim_{n\to\infty}(\epsilon_n/\epsilon_{n-i})=0\quad\text{for all $i\geq
1$,}$$ and

$$  \epsilon_n/\epsilon_{n-i}=o(\epsilon_n/\epsilon_{n-j})\quad\text{as
$n\to\infty$,\quad for $j<i$.}$$
 Consequently, expanding in \eqref{eq:123} the product $\prod^k_{i=1}(\epsilon_n-\epsilon_{n-i})$, we have

\begin{align}\label{eq:200}\prod^k_{i=1}(\epsilon_n-\epsilon_{n-i})&=\prod^k_{i=1}\bigg(-\epsilon_{n-i}
[1-\epsilon_n/\epsilon_{n-i}]\bigg)
 \nonumber \\
&= (-1)^k
\bigg(\prod^k_{i=1}\epsilon_{n-i}\bigg)[1+O(\epsilon_n/\epsilon_{n-1})]\quad\text{as
$n\to\infty$.}\end{align}
 Substituting \eqref{eq:200}
in \eqref{eq:123},  and defining

\begin{equation}\label{eq:271}
D_n=\frac{\widehat{D}_n}{p'_{n,k}(z_n)}, \quad
E_n=\frac{\widehat{E}_n}{p'_{n,k}(z_n)},\end{equation}  we obtain

\begin{equation}\label{eq:231} \epsilon_{n+1}=(-1)^k
D_n\bigg(\prod^k_{i=0}\epsilon_{n-i}\bigg)[1+O(\epsilon_n/\epsilon_{n-1})]+E_n\epsilon_n^2\quad
\text{as $n\to\infty$.}
\end{equation}
Dividing both sides of \eqref{eq:231} by
$\prod^k_{i=0}\epsilon_{n-i}$, and defining

\begin{equation}\label{eq:233}
\sigma_n=\frac{\epsilon_{n+1}}{\prod^k_{i=0}\epsilon_{n-i}},
\end{equation}
 we have

\begin{equation}\label{eq:235}
 \sigma_n=(-1)^kD_n[1+O(\epsilon_n/\epsilon_{n-1})] +E_n \sigma_{n-1}\epsilon_{n-k-1}\quad \text{as
$n\to\infty$.}
\end{equation}

\sloppypar
Now, by  \eqref{eqlim}, \eqref{eq123}, and \eqref{ppff},
 \begin{equation}\label{eq:237}
\lim_{n\to\infty}D_n=-\frac{1}{(k+1)!}\frac{f^{(k+1)}(\alpha)}{f'(\alpha)},\quad
\lim_{n\to\infty}E_n=\frac{f^{(2)}(\alpha)}{2f'(\alpha)}.\end{equation}
Because  $\lim_{n\to\infty}D_n$ and $\lim_{n\to\infty}E_n$ are
finite, and because $\lim_{n\to\infty}(\epsilon_n/\epsilon_{n-1})=0$ and
$\lim_{n\to\infty}\epsilon_{n-k-1}=0$,  it follows that there
exist a positive integer $N$ and positive constants $\beta<1$ and
$D$,  with $|E_n\epsilon_{n-k-1}|\leq \beta$ when $n> N$, for
which \eqref{eq:235} gives

\begin{equation}\label{eq:241}
 |\sigma_n|\leq D+\beta|\sigma_{n-1}|\quad\text{for all
$n> N$.}\end{equation} Using \eqref{eq:241}, it is easy to show
that
$$|\sigma_{N+s}|\leq
D\frac{1-\beta^s}{1-\beta}+\beta^s|\sigma_N|,\quad s=1,2,\ldots,$$
which, by the fact that $\beta<1$, implies that $\{\sigma_n\}$ is
a bounded sequence. Making use of this fact, we
 have $\lim_{n\to\infty}E_n \sigma_{n-1}\epsilon_{n-k-1}=0$. Substituting this
 in \eqref{eq:235}, and  invoking \eqref{eq:237}, we next obtain
$\lim_{n\to\infty}\sigma_n=(-1)^k\lim_{n\to\infty}D_n=L$, which is
precisely \eqref{eq:130}.

That the order of the method is $s_k$, as defined in the statement
of the theorem, follows from Traub
 \cite[Chapter 3]{Traub:1964:IMS}.  A  weaker version can be proved by
letting $\sigma_n=L$ for all $n$ and showing that
$|\epsilon_{n+1}|=Q|\epsilon_n|^{s_k}$ is possible for $s_k$ a
solution to the equation $s^{k+1}=\sum^k_{i=0}s^i$ and
$Q=|L|^{(s_k-1)/k}$. The  proof of this is easy and is left  to
the reader. This completes the proof of part 1 of the theorem.

When $f(z)$ is a polynomial of degree at most $k$, we first observe that $f^{(k+1)}(z)=0$ for all $z$, which implies that $p_{n,k}(z)= f(z)$ for all $z$, hence also
$p'_{n,k}(z)= f'(z)$ for all $z$.  Therefore, we have that
$p'_{n,k}(z_n)=f'(z_n)$ in the recursion  of \eqref{eq:gs}. Consequently,  \eqref{eq:gs} becomes
$$ z_{n+1}=z_n-\frac{f(z_n)}{f'(z_n)}, \quad n=k,k+1,\ldots, $$
which is the recursion for the Newton--Raphson method.  Thus,
\eqref{eq:177} follows. This completes the  proof of part 2 of the theorem.
\hfill $\blacksquare$

\section{Numerical examples}\label{se3}
\setcounter{theorem}{0} \setcounter{equation}{0}

In this section, we present two numerical examples that we treated with our method.
 Our computations were done in
quadruple-precision arithmetic (approximately 35-decimal-digit
accuracy).  Note that in
order to verify the theoretical results concerning iterative
methods of order greater that unity, we need to use computer
arithmetic of high precision (preferably, of variable precision,
if available) because the number of correct significant decimal
digits in the $z_n$ increases dramatically from one iteration to the next as we
are approaching the solution.

In both examples below, we take $k=2$. We choose $z_0$ and $z_1$  and compute $z_2$
using one step of the secant method, namely,
\begin{equation} \label{eq:5-1}
z_2=z_1-\frac{f(z_1)}{f[z_0,z_1]}.\end{equation}
 Following that,
we compute $z_3,z_4,\ldots,$ via
\begin{equation} \label{eq:5-2}
 z_{n+1}=z_n-\frac{f(z_n)}
{f[z_n,z_{n-1}]+f[z_n,z_{n-1},z_{n-2}](z_n-z_{n-1})},\quad
n=2,3,\ldots\ .\end{equation}
In our  examples, we have carried out our computations for several sets of $z_0,z_1$,  and we have observed essentially the same behavior that we observe in Tables \ref{table10} and \ref{table20}.

\begin{example}{\em
 Consider  $f(z)=0$, where
$f(z)=z^3-8$, whose solutions are  $\alpha_r=2e^{\mrm{i}2\pi r/3}$, $r=0,1,2$.
We  would like to obtain
the root $\alpha_1=2e^{\mrm{i}2\pi/3}=-1+\mrm{i}\sqrt{3}$.
   We chose  $z_0=2\mrm{i}$ and  $z_1=-2+2\mrm{i}$.
The results of our computations are given in Table \ref{table10}.

\begin{table}[tb]
\begin{center}
$$
\begin{array}{||c||c|c|c|c||}
\hline
 n&|\epsilon_n|&\displaystyle
 \frac{\epsilon_{n+1}}{\epsilon_n\epsilon_{n-1}\epsilon_{n-2}}
 & \displaystyle\frac{\log|\epsilon_{n+1}/\epsilon_n|}{\log|\epsilon_{n}/\epsilon_{n-1}|}\\
\hline\hline
  0&      1.035D+00&      -   &   -   \\
  1&      1.035D+00&      -  &   -   \\
  2&      4.808D-01&     -8.972D-02   +\mrm{i}\    1.015D-01~~~&      2.516\\
  3&      6.979D-02&      1.224D-01      -\mrm{i}~2.727D-02&      1.437\\
  4&      4.355D-03&      1.009D-01      -\mrm{i}~4.079D-02&      2.023\\
  5&      1.591D-05&      4.561D-02      -\mrm{i}~9.794D-02&      1.839\\
  6&      5.223D-10&      3.793D-02      -\mrm{i}~7.268D-02&      1.839\\
  7&      2.967D-18&      3.741D-02      -\mrm{i}~7.579D-02&      1.838\\
  8&      2.083D-33&      *    *&       **\\
  9&      0.000D+00&           *      *     &           **\\
   \hline

\end{array}
$$

\caption{\label{table10} Results obtained by applying the
generalized secant method with $k=2$, as shown in \eqref{eq:5-1}
and \eqref{eq:5-2}, to the equation $z^3-8=0$, to compute the root $\alpha_1=-1+\mrm{i}\sqrt{3}$.}

\end{center}
\end{table}

From \eqref{eq:130} and \eqref{eq:135} in Theorem \ref{th1}, we should have

$$
\lim_{n\to\infty}\frac{\epsilon_{n+1}}{\epsilon_{n}\epsilon_{n-1}\epsilon_{n-2}}=
\frac{(-1)^{3}}{3!}\,\frac{f'''(\alpha_1)}{f'(\alpha_1)}=\frac{1}{24}(1-\mrm{i}\sqrt{3})
=0.04166\cdots-\mrm{i}\,0.07216\cdots$$
and $$ \lim_{n\to\infty}\frac{\log|\epsilon_{n+1}/\epsilon_{n}|}
{\log|\epsilon_{n}/\epsilon_{n-1}|}=s_2=1.83928\cdots,$$ and these
seem to be confirmed in Table \ref{table10}. Also, in infinite-precision arithmetic, $x_9$
 should have close to 60 correct significant
figures; we do not see this in  Table \ref{table10}
due to the fact that the arithmetic we have used to generate Table
\ref{table10} can provide an accuracy of at most 35 digits.
}
\end{example}

\begin{example}{\em
 Consider  $f(z)=0$, where
$f(z)=\sin(\mrm{i}z)-\cos z$. $f(z)$  has infinitely many roots
$\alpha_r=(1-\mrm{i})(\pi/4+r\pi)$, $r=0,\pm1,\pm2,\ldots$.
We would like to obtain the root $\alpha_0=(1-\mrm{i})\pi/4$.
We chose $z_0=1.5-\mrm{i}1.3$ and  $z_1=0.6-\mrm{i}0.5$.
The results of our computations are given in Table \ref{table20}.
\begin{table}[tb]
\begin{center}
$$
\begin{array}{||c||c|c|c|c||}
\hline
 n&|\epsilon_n|&\displaystyle
 \frac{\epsilon_{n+1}}{\epsilon_n\epsilon_{n-1}\epsilon_{n-2}}
 & \displaystyle\frac{\log|\epsilon_{n+1}/\epsilon_n|}{\log|\epsilon_{n}/\epsilon_{n-1}|}\\
\hline\hline
 0&      6.608D-01&      -    &      -\\
  1&      3.403D-01&      -&      -\\
  2&      1.341D-01&      3.163D-01    +\mrm{i}\,  1.397D-01&      2.743\\
  3&      1.043D-02&      1.466D-01     -\mrm{i}\,1.846D-01&      1.774\\
  4&      1.122D-04&     -2.943D-03     -\mrm{i}\,1.117D-01~~~&      1.934\\
  5&      1.755D-08&      9.223D-03     -\mrm{i}\,1.614D-01&      1.766\\
  6&      3.320D-15&     -7.686D-04     -\mrm{i}\,1.658D-01~~~&      1.857\\
  7&      1.084D-27&      **      &      ** \\
  8&      9.630D-35&      **       &     **\\
    \hline
\end{array}
$$

\caption{\label{table20} Results obtained by applying the
generalized secant method with $k=2$, as shown in \eqref{eq:5-1}
and \eqref{eq:5-2}, to the equation $\sin(\mrm{i}z)-\cos z=0$, to compute the root $\alpha_0=(1-\mrm{i})\pi/4$.}

\end{center}
\end{table}

From \eqref{eq:130} and \eqref{eq:135} in Theorem \ref{th1}, we should have

$$
\lim_{n\to\infty}\frac{\epsilon_{n+1}}{\epsilon_{n}\epsilon_{n-1}\epsilon_{n-2}}=
\frac{(-1)^{3}}{3!}\,\frac{f'''(\alpha_1)}{f'(\alpha_1)}=-\frac{\mrm{i}}{6}=
-\mrm{i}\,0.16666\cdots$$
and
$$ \lim_{n\to\infty}\frac{\log|\epsilon_{n+1}/\epsilon_{n}|}
{\log|\epsilon_{n}/\epsilon_{n-1}|}=s_2=1.83928\cdots,$$ and these
seem to be confirmed in Table \ref{table20}. Also, in infinite-precision arithmetic, $x_8$
 should have close to 50 correct significant
figures; we do not see this in  Table \ref{table20}
due to the fact that the arithmetic we have used to generate Table
\ref{table20} can provide an accuracy of at most 35 digits.
}
\end{example}

\noindent{\bf Remark.} Before closing, we would like to discuss the issue of estimating the relative errors   $|\epsilon_n/\alpha|$ in  the $z_n$. This should help the reader when studying the numerical results included
in Tables \ref{table10} and \ref{table20}.
Starting with \eqref{eq:130} and   \eqref{eq:135},  we first note that, for all large $n$,
$$ |\epsilon_{n+1}|\approx |L|^{(s_k-1)/k}|\epsilon_n|^{s_k}.$$
Therefore,  assuming also that $\alpha\neq0$, we have
$$ |\epsilon_{n+1}/\alpha|\approx  D |\epsilon_{n}/\alpha|^{s_k},\quad D=\big(|L|^{1/k} |\alpha|\big)^{s_k-1}.$$
 Now, if $z_n$ has $q>0$ correct significant figures, we have $|\epsilon_{n}/\alpha|=O(10^{-q})$. If, in addition,  $D=O(10^r)$ for some $r$, then we will have

$$|\epsilon_{n+1}/\alpha|\approx O(10^{r-qs_k}).$$
 For simplicity, let us consider the case $r=0$, which is practically what we have in the two examples we have treated. Then
 $z_{n+1}$ has  approximately   $qs_k$ correct significant decimal digits. That is, if $z_n$ has $q$ correct significant decimal digits, then, due to the fact that $s_k>1$, $z_{n+1}$ will have $s_k$ times as many  correct significant decimal digits as~$z_n$.


\end{document}